\documentclass[12pt]{article}
\usepackage{amssymb}
\usepackage{latexsym}
\usepackage{amsfonts}
\usepackage{amsmath}
\newtheorem{thrm}{Theorem}[section]
\newtheorem{lemma}[thrm]{Lemma}
\newtheorem{prop}[thrm]{Proposition}

\newtheorem{remark}[thrm]{Remark}
\newtheorem{assumption}{Assumption}

\numberwithin{equation}{section}

\usepackage[dvips]{color}

\def\P{\mathbb{P} }
\def\R{\mathbb{R} }

\def\la{\langle}
\def\ra{\rangle}

\begin{document}
\allowdisplaybreaks

\allowdisplaybreaks
\begin{doublespace}
\title{\Large\bf On properties of a class of strong limits for
supercritical superprocesses}
\author{ \bf  Yan-Xia Ren\footnote{The research of this author is supported by NSFC (Grant No.   11671017 and 11731009).\hspace{1mm} } \hspace{1mm}\hspace{1mm}
Renming Song\thanks{Research supported in part by a grant from the Simons
Foundation (\#429343, Renming Song).} \hspace{1mm}\hspace{1mm} and \hspace{1mm}\hspace{1mm}
Rui Zhang\footnote{Research supported in part by NSFC (Grant No. 11601354). Corresponding author.}
\hspace{1mm} }
\date{}
\maketitle

\begin{abstract}
Suppose that $X=\{X_t, t\ge 0; \P_{\mu}\}$ is a supercritical superprocess in
a locally compact separable metric space $E$.  Let $\phi_0$ be a positive
 eigenfunction corresponding to the first eigenvalue $\lambda_0$ of
the generator of the mean semigroup of $X$. Then
$M_t:=e^{-\lambda_0t}\langle\phi_0, X_t\rangle$ is a positive martingale.
Let $M_\infty$ be the limit of $M_t$. It is known that
$M_\infty$ is non-degenerate iff  the $L\log L$ condition is  satisfied.
When the $L\log L$ condition may not be satisfied, we recently proved in (arXiv:1708.04422) that
there exist a non-negative
function $\gamma_t$ on $[0, \infty)$
and  a non-degenerate random variable
$W$ such that for any
finite nonzero Borel measure $\mu$ on $E$,
$$
\lim_{t\to\infty}\gamma_t\langle \phi_0,X_t\rangle
=W,\qquad\mbox{a.s.-}\P_{\mu}.
$$
In this paper, we mainly investigate properties of $W$. We prove that $W$ has strictly positive density on $(0,\infty)$.
We also investigate the small value probability and tail probability problems of $W$.
\end{abstract}

\medskip


\medskip


\bigskip

\baselineskip=6.0mm
\section{Background and our model}

Consider a supercritical Galton-Watson process $\{Z_n, n\ge 0\}$
with offspring distribution $\{p_n: n\ge 0\}$. In 1968, Seneta \cite{Sen68}
proved that there exists
a sequence of positive numbers $\{c_n,n\ge1\}$ such that
$c_nZ_n$ converges in distribution to a non-degenerate random variable $W$; then
Heyde \cite{{Heyde70}} strengthened this convergence to almost sure convergence.
Since then, the problem of finding $\{c_n,n\ge1\}$ such that $c_nZ_n$ converges
to a non-degenerate limit is called the Seneta-Heyde norming problem, $\{c_n,n\ge1\}$
are called the norming constants. Harris \cite{Harris} proved that, when $\{p_n: n\ge 0\}$
has a finite second moment, the distribution of $W$, restricted to $(0, \infty)$,
is absolutely continuous; then  Stigum \cite{stigum}
extended this to the case when $\{p_n: n\ge 0\}$
satisfies the $L\log L$ condition. Finally, Athreya \cite{Athreya71} proved that the same conclusion holds for all
supercritical Galton-Watson processes.
As for other properties of $W$,
\cite{MO} discussed the small value probability problem for $W$, i.e., the rate at which
the probability $P(0<W\le r)$ tends to 0 as $r\to 0$; \cite{BB} studied the tail probability problem for $W$, i.e.,
  the rate at which $P(W> r)$ tends to 0 as $r\to\infty$,
 under the assumption that there
exists $N>0$ such that $p_n=0$ for all $n\ge N$.

For supercritical multitype Galton-Watson process, Jones \cite{Jones} studied the corresponding
small value probability problem and tail probability problem.
Hering \cite{Hering} established the corresponding results  for supercritical branching Markov processes.
In the recent paper \cite{RSZ6},
we studied the Seneta-Heyde type limit problem for supercritical superprocesses:
Suppose $\{X_t, t\ge 0; \P_{\mu}\}$ is a supercritical superprocess on $E$, we proved that, under certain conditions,
there exist a non-negative function $\gamma_t$ on $[0, \infty)$
and a non-degenerate
random variable $W$ such that for all finite Borel measure $\mu$ on $E$,
\begin{equation*}
\lim_{t\to\infty}\gamma_t\langle \phi_0,X_t\rangle
=W,\qquad\mbox{a.s.-}\P_{\mu},
\end{equation*}
where $\phi_0$ is a positive eigenfunction of the infinitesimal generator of the mean semigroup
of $X$ corresponding to the first eigenvalue $\lambda_0$. The main goals of this paper are to
study the absolute continuity of $W$, when restricted to $(0, \infty)$, and to study the
small value probability problem and tail probability problem for $W$.

\subsection{Superprocesses}

We first introduce the setup of this paper.
Suppose that $E$ is a locally compact separable metric space,
and $\partial$ is a separate point not contained in $E$.
We will use $E_{\partial}$ to denote $E\cup\{\partial\}$.
Suppose that $m$ is a
$\sigma$-finite Borel measure on $E$ with full support.
We will use
$\mathcal{B}(E)$ ($\mathcal{B}^+(E)$) to denote the
family of (non-negative)  Borel  functions on $E$,
$\mathcal{B}_b(E)$ ($\mathcal{B}_b^+(E)$) to denote the
family of (non-negative) bounded Borel  functions on $E$,
and $C(E)$ to denote the family of continuous functions on $E$.
We assume that $\xi=\{\Omega^0, \mathcal{H}, \mathcal{H}_t, \xi_t,\Pi_x, \zeta\}$  is a
Hunt process on $E$,
where $\{\mathcal{H}_t: t\ge 0\}$ is the minimal filtration of $\xi$
satisfying the usual conditions and $\zeta:=
 \inf\{t>0: \xi_t=\partial\}$ is the lifetime of $\xi$.
We will use $\{P_t:t\geq 0\}$ to denote the semigroup of $\xi$.

The superprocess $X=\{X_t:t\ge 0\}$ we are going to work with is determined by two parameters:
a spatial motion $\xi=\{\xi_t, \Pi_x\}$ on $E$
which is a Hunt process,
and a branching mechanism $\varphi$ of the form
\begin{equation}\label{e:branm}
\varphi(x,s)=-\alpha(x)s+\beta(x)s^2+\int_{(0,+\infty)}(e^{-s r}-1+s r)n(x,{\rm d}r),
\quad x\in E, \, s\ge 0,
\end{equation}
where $\alpha\in \mathcal{B}_b(E)$, $\beta\in \mathcal{B}_b^+(E)$ and $n$ is a kernel from $E$ to $(0,\infty)$ satisfying
\begin{equation}\label{n:condition}
  \sup_{x\in E}\int_{(0,+\infty)} (r\wedge r^2) n(x,{\rm d}r)<\infty.
\end{equation}
It follows from the above assumptions that  there exists $M>0$ such that
\begin{equation}\label{e:M}
|\alpha(x)|+\beta(x)+
\int_{(0,+\infty)} (r\wedge r^2) n(x,{\rm d}r)\le M.
\end{equation}

Let ${\cal M}_F(E)$
be the space of finite measures on $E$, equipped with the topology of weak convergence.
The superprocess $X$  with spatial motion $\xi$ and branching mechanism $\varphi$
is a  Markov process taking values in ${\cal M}_F(E)$.
The existence of such superprocesses is well-known,
see \cite{Li11} or \cite{E.B.}, for instance.
For any $\mu \in \mathcal{M}_F(E)$, we denote the
law
of $X$ with initial configuration $\mu$ by $\P_\mu$.
As usual,
$\langle f,\mu\rangle:=\int_E f(x)\mu(dx)$
and $\|\mu\|:=\langle 1,\mu\rangle$.
Throughout this paper, a real-valued function $u(t,x)$ on $[0,\infty)\times E_\partial$
is said to be locally bounded if, for any $t>0$, $\sup_{s\in[0,t], x\in E_\partial}|u(s, x)|<\infty$.
Any function $f$ on $E$ is automatically extended to $E_\partial$ by setting $f(\partial)=0$.
According to \cite[Theorem 5.12]{Li11}, there is a Hunt process
$X=\{\Omega, {\cal G}, {\cal G}_t, X_t, \P_\mu\}$ taking values in  $\mathcal{M}_F(E)$
such that for every
$f\in \mathcal{B}^+_b(E)$ and $\mu \in \mathcal{M}_F(E)$,
\begin{equation}
  -\log \P_\mu\left(e^{-\langle f,X_t\rangle}\right)=\langle V_tf,\mu\rangle,
\end{equation}
where $V_tf(x)$ is the unique locally bounded non-negative solution to the equation
\begin{equation}\label{V_t}
   V_tf(x)+\Pi_x\int_0^t\varphi(\xi_s, V_{t-s}f(\xi_s)){\rm d}s=\Pi_x f(\xi_t), \quad x\in E_\partial,
\end{equation}
where we use the convention that $\varphi(\partial,r)=0$ for all $r\ge0$.
Since $f(\partial)=0$, we have $V_tf(\partial)=0$ for any $t\ge 0$.
In this paper, the superprocess we deal with is always this Hunt realization.

For any $f\in\mathcal{B}_b(E)$ and $(t, x)\in (0, \infty)\times E$, we define
\begin{equation}\label{61.26}
   T_tf(x):=\Pi_x \left[e^{\int_0^t\alpha(\xi_s){\rm d}s}f(\xi_t)\right].
\end{equation}
It is well known that $T_tf(x)=\P_{\delta_x}\langle f,X_t\rangle$ for every $x\in E$.

We will always assume that there exists a family of
continuous and strictly positive functions $\{p(t,x,y):t>0\}$ on $E\times E$ such that for any $t>0$ and non-negative function $f$ on $E$,
$ P_tf(x)=\int_E p(t,x,y)f(y)m({\rm d}y)$.
Define
$$
a_t(x):=\int_E p(t,x,y)^2\,m({\rm d}y),\qquad \hat{a}_t(x):=\int_E p(t,y,x)^2\,m({\rm d}y).$$
Our first assumption is
\begin{assumption}\label{assum1'}
\begin{description}
  \item[(i)]  For any $t>0$, $\int_E p(t,x,y)\,m({\rm d}x)\le 1$.
  \item[(ii)] For any $t>0$, we have
  \begin{equation}\label{61.17}
     \int_E a_t(x)\,m({\rm d}x)=\int_E\hat{a}_t(x)\,
m({\rm d}x)=\int_E\int_E p(t,x,y)^2\,m({\rm d}y)\,m({\rm d}x)<\infty.
  \end{equation}
  Moreover, the functions $x\to a_t(x)$ and $x\to\hat{a}_t(x)$ are continuous on $E$.
\end{description}
\end{assumption}

Note that, in Assumption \ref{assum1'}(i), the integration is with respect to the first
space variable. This implies that the dual semigroup $\{{\widehat P}_t:t\ge 0\}$
of $\{P_t:t\ge 0\}$ with respect to $m$
is sub-Markovian.
By H\"older's inequality, we have
\begin{equation}\label{61.1}
 p(t+s,x,y)=\int_E p(t,x,z)p(s,z,y)\,m({\rm d}z)\le (a_t(x))^{1/2}(\hat{a}_s(y))^{1/2}.
\end{equation}
$\{P_t:t\ge 0\}$
and $\{\widehat{P}_t:t\ge 0\}$ are both strongly continuous
contraction semigroups on $L^2(E, m)$,
see \cite{RSZ3} for a proof.
We will use $\langle\cdot, \cdot\rangle_m$ to denote the inner product in $L^2(E, m)$.
Since $p(t, x, y)$ is continuous in $(x, y)$,
by \eqref{61.1}, Assumption \ref{assum1'}(ii) and
the dominated convergence theorem, we have that,
for any $f\in L^2(E,m)$,  $P_tf$  and $\widehat{P}_tf$ are continuous.

It follows from Assumption \ref{assum1'}(ii) that,
for each $t>0$, $\{P_t\}$ and $\{\widehat{P}_t\}$ are compact operators on $L^2(E,m)$.
Let $\widetilde{L}$ and $\widehat{\widetilde{L}}$ be the infinitesimal generators of the semigroups $\{P_t\}$ and $\{\widehat{P}_t\}$ in $L^2(E,m)$ respectively.
Define $\widetilde{\lambda}_0:=\sup \Re(\sigma(\widetilde{L}))=
\sup\Re(\sigma(\widehat{\widetilde{L}}))$, where $\Re$ stand for the real part of a complex number.
By Jentzsch's theorem (\cite[Theorem V.6.6]{Sch}),
$\widetilde{\lambda}_0$ is an eigenvalue of multiplicity 1 for both $\widetilde{L}$ and
$\widehat{\widetilde{L}}$. Let $\widetilde{\phi}_0$  and $\widetilde{\psi}_0$
be, respectively,
eigenfunctions of $\widetilde{L}$ and $\widehat{\widetilde{L}}$
corresponding to $\widetilde{\lambda}_0$.
$\widetilde{\phi}_0$  and $\widetilde{\psi}_0$ can be chosen
be strictly positive $m$-almost
everywhere with $\|\widetilde{\phi}_0\|_2=1$ and $\langle \widetilde{\phi}_0, \widetilde{\psi}_0\rangle_m=1$.
Thus for $m$-almost every $x\in E$,
$$
e^{\widetilde{\lambda}_0}\widetilde{\phi}_0(x)=P_1\widetilde{\phi}_0(x),
\qquad
e^{\widetilde{\lambda}_0}\widetilde{\psi}_0(x)=\widehat{P}_1\widetilde{\psi}_0(x).
$$
Hence, by the continuity of $P_1\widetilde{\phi}_0$ and
$\widehat{P}_1\widetilde{\psi}_0$,
 $\widetilde{\phi}_0$ and $\widetilde{\psi}_0$ can be chosen to be continuous and strictly
positive everywhere on $E$.

Our second assumption is that $\{P_t:t\ge 0\}$
and $\{\widehat{P}_t:t\ge 0\}$ are intrinsically ultracontractive.

\begin{assumption}\label{assum63}
\begin{description}
 \item[(i)] $\widetilde{\phi}_0$ is bounded.

 \item[(ii)]
The semigroups $\{P_t,t\ge0\}$ and $\{\widehat{P}_t:t\ge 0\}$ are intrinsically ultracontractive,
that is, there exists $c_t>0$ such that
\begin{equation}\label{newcondition2}
  p(t,x,y)\le c_t\widetilde{\phi}_0(x)\widetilde{\psi}_0(y).
\end{equation}
\end{description}
\end{assumption}

In \cite{RSZ3}, we have given many examples of Hunt processes satisfying Assumptions
\ref{assum1'}--\ref{assum63}. For example, if $E$ be a bounded Lipschitz domain,
and $\xi$ is the subprocess of a diffusion process whose generator is a uniformly elliptic
second order differential operator, then $\xi$ satisfies
Assumptions \ref{assum1'}--\ref{assum63}, see \cite{DS}.

By using the boundedness of $\alpha$ and assumptions on $\xi$, we have proved in \cite[Lemma 2.1]{RSZ4} that
the semigroup $\{T_t\}$ has a continuous and strictly positive density $q(t, x, y)$ with
respect to the measure $m$, that is, for any $f\in \mathcal{B}_b(E)$,
$$
  T_tf(x)=\int_E q(t,x,y)f(y)\,m({\rm d}y).
$$
and, for any $t>0$, $q(t, x, y)$ is continuous in $(x, y)$, and
\begin{equation}\label{comp0}
e^{-Mt}p(t,x,y) \le q(t,x,y)\le e^{Mt}p(t,x,y), \quad (t, x, y)\in (0, \infty)\times E\times E,
\end{equation}

Let $\{\widehat{T}_t,t>0\}$ be the dual semigroup of $\{T_t,t>0\}$ in $L^2(E,m)$, that is,
for any $f,g\in L^2(E,m)$,
\begin{equation*}
\widehat{T}_tf(x)=\int_E q(t,y,x)f(y)\,m(dy).
\end{equation*}
It follows from Assumption \ref{assum1'}(ii) and \eqref{comp0} that

\begin{equation*}
\int_E\int_E q^2(t,x,y)\,m(x)\,m(dy)\le e^{2Mt}\int_E\int_E p^2(t,x,y)\,m(x)\,m(dy)<\infty.
\end{equation*}
Thus using the same analysis as that used before Assumption \ref{assum63} we can get the following conclusion: for any $t>0$, $T_t$ and $\widehat{T}_t$ are compact operators
on $L^2(E,m)$. Let $L$ and $\widehat{L}$ be the infinitesimal generators of $\{T_t\}$
and $\{\widehat{T}_t\}$ in $L^2(E,m)$ respectively. Define
$\lambda_0:=\sup \Re(\sigma(L))=\sup\Re(\sigma(\widehat{L}))$.
$\lambda_0$ is an eigenvalue of multiplicity one for both
$L$ and $\widehat{L}$. Let $\phi_0$ and $\psi_0$ be, respectively, eigenfunctions of $L$ and $\widehat{L}$ corresponding to $\lambda_0$. $\phi_0$ and $\psi_0$ can be chosen to be continuous and strictly positive everywhere with $\|\phi_0\|_2=1$, $\langle \phi_0,\psi_0\rangle_m=1$.

Using Assumption \ref{assum63}, the boundedness of $\alpha$ and an
argument similar to that used in the proof of \cite[Theorem 3.4]{DS},
one can show the following:
\begin{description}
\item[(i)] $\phi_0$ is bounded.

\item[(ii)]
The semigroup $\{T_t,t\ge0\}$ and $\{\widehat{T}_t,t>0\}$ are intrinsically ultracontractive,
that is, there exists $c_t>0$ such that
\begin{equation}\label{6condition2}
  q(t,x,y)\le c_t\phi_0(x)\psi_0(y).
\end{equation}
\end{description}

Define $q_t(x):=\P_{\delta_x}(\|X_t\|=0)$ and  $\mathcal{E}:=\{\exists t>0, \|X_t\|=0\}$.
Note that $q_t(x)$ is non-decreasing in $t$. Thus the limit
\begin{equation*}
q(x):=\lim_{t\to\infty}q_t(x)=\P_{\delta_x}(\mathcal{E})
\end{equation*}
exist. We call $q(x)$ the extinction probability of the superprocess.  Let $v(x):=-\log q(x)$.
It follows from the branching property of $X$ that $\P_{\mu}(\mathcal{E})=e^{-\la v,\mu\ra}$.
The main interest of this paper is on supercritical superprocesses, so we assume that
\begin{assumption}\label{assum2}
  $\lambda_0>0$.
\end{assumption}

We also assume that
\begin{assumption}\label{assum4}
There exists $t_0>0$ such that
\begin{equation}\label{4.2'}
  \inf_{x\in E} q_{t_0}(x)>0.
\end{equation}
\end{assumption}

Assumption \ref{assum4} guarantees that $\|v\|_\infty\le \sup_{x\in E}(-\log q_{t_0}(x))
<\infty$, thus $v$ is a bounded function.
In \cite[Section 2.2]{RSZ4}, we gave
a sufficient condition  for Assumption \ref{assum4}. In particular, if $\inf_{x\in E}\beta(x)>0$, then Assumption \ref{assum4} holds.
Under Assumptions \ref{assum1'}--\ref{assum4}, we have proven in
\cite[Lemma 3.1]{RSZ6} that
$q(x)<1,$ for all $x\in E$,
which is a reflection of supercriticality.

\subsection{Main results}

Define
\begin{equation*}
M_t:=e^{-\lambda_{0}t}\langle \phi_0,X_{t}\rangle,\quad t\ge 0.
\end{equation*}
It follows from the Markov property that, for every $\mu\in\mathcal{M}_F(E)$,
$\{M_t, t\ge 0\}$
is a non-negative $\mathbb{P}_{\mu}$-martingale.
Thus $\{M_t, t\ge 0\}$ has a $\mathbb{P}_{\mu}$-a.s. finite limit denoted as $M_\infty$.
According to \cite{LRS09}, $M_\infty$ is non-degenerate if and only if the $L\log L$ condition holds. When $M_\infty$ is a non-degenerate random variable,
$X_t$ grows exponentially and the growth rate is
$e^{\lambda_0 t}$. When $M_\infty$ is a degenerate random variable, $e^{\lambda_0 t}$
is no longer the growth rate of $X_t$. In \cite{RSZ6}, we proved
that, when the $L\log L$ condition may not be satisfied, the growth rate of $X_t$
is $e^{\lambda_0 t}L(t)$,
where $L(t)$ is a slowly varying function after some transform.
Now we state the main results of \cite{RSZ6}.

In \cite{RSZ6}, we proved that there exists a family of non-negative functions
$\{\eta_t(x),t\ge 0\}$, satisfying $0\le \eta_t(x)\le v(x)$, such that
\begin{equation*}
\eta_t(x)=V_s(\eta_{t+s})(x),\qquad t, s\ge0, x\in E.
\end{equation*}
Furthermore, $\eta_0$ is not identically 0, is also not identically equal to $v$. Let
$\gamma_t=\langle \eta_t, \psi_0\rangle_m$, then
\begin{equation*}
\lim_{t\to\infty}\frac{\gamma_{t}}{\gamma_{t+s}}=e^{\lambda_0s},\quad\forall s\ge 0,\end{equation*}
and the following assertions are valid.

\begin{thrm}\emph{{\bf\cite[Theorem 1.2]{RSZ6}}}.
There exists a non-degenerate random variable $W$ such that for all $\mu\in\mathcal{M}_F(E)$,
\begin{equation*}
\lim_{t\to\infty}\gamma_t\langle \phi_0,X_t\rangle
=W,\qquad
\mbox{a.s.-}\P_{\mu}
\end{equation*}
and
\begin{equation*}
\P_{\mu}(W=0)=e^{-\la v,\mu\ra},\qquad \P_{\mu}(W<\infty)=1.
\end{equation*}
\end{thrm}

Define a new measure $n^{\phi_0}(x,{\rm d}r)$ by
\begin{equation*}
\int_0^\infty f(r)n^{\phi_0}(x,{\rm d}r)=\int_0^\infty f(r\phi_0(x))n(x,{\rm d}r).
\end{equation*}
Let $l(x):=\int_1^\infty r\ln r \,n^{\phi_0}(x, {\rm d}r)$. Necessary and sufficient
conditions for $M_\infty$ to be non-degenerate are as follows:

\begin{thrm}\emph{\bf {\cite[Theorem1.3]{RSZ6}}}.
The following are equivalent:
\begin{description}
  \item[(1)] for some $\mu\in{\cal M}_F(E)$, $M_\infty$ is a non-degenerate random variable under $\mathbb{P}_{\mu}$;
  \item[(2)]for every non-zero $\mu\in{\cal M}_F(E)$, $M_\infty$ is a non-degenerate random variable under $\mathbb{P}_{\mu}$;
  \item[(3)] $l_0:=\lim_{t\to\infty}e^{\lambda_0t}\gamma_t<\infty$;
  \item[(4)] (\emph{\bf {$L\log L$ criterion:}}) $\int_E\psi_0(x)l(x)m({\rm d}x)<\infty;$
  \item[(5)] for some non-zero $\mu\in{\cal M}_F(E)$, $\P_{\mu}W<\infty$;
  \item[(6)] for every non-zero $\mu\in{\cal M}_F(E)$, $\P_{\mu}W<\infty$.
\end{description}
\end{thrm}

It follows from \cite[Remark 1.1]{RSZ6} that, $\int_E\psi_0(x)l(x)m({\rm d}x)<\infty$ if
and only if
\begin{equation}\label{llog2'}
\int_E\phi_0(x)\psi_0(x)m({\rm d}x)\int^\infty_1(r\ln r)n(x, {\rm d}r)<\infty.
\end{equation}

The main purpose of this paper is to further study properties of $W$: including whether
$W$ has a density function, the small value probability problem and the tail problem
for $W$. The main results of this paper are as follows.

\begin{thrm}\label{main-dens}
For any non-zero $\mu\in\mathcal{M}_F(E)$, under $\P_{\mu}$,
the restriction of the random variable $W$ on $(0,\infty)$ has a strictly positive
density.
\end{thrm}

In Subsection \ref{sec:est}, we will introduce another semigroup $\{T_t^*\}$ with largest
eigenvalue $\lambda_0^*<0$. Define
\begin{equation*}
\epsilon_0:=\frac{-\lambda_0^*}{\lambda_0}.
\end{equation*}
Let
\begin{equation}\label{L(t)}
L(t)=e^{-\lambda_0t}\gamma_t,\end{equation}
then
\begin{equation*}\lim_{t\to\infty}\frac{L(t+s)}{L(t)}=1.
\end{equation*}
Define
\begin{equation}
\label{tilde-L}\widetilde{L}(\theta):=L(\log \theta/\lambda_0),\quad\theta\ge 1,
\end{equation}
then $\widetilde{L}$ is a slowly varying function at $\infty$.

\begin{thrm}\label{main-tail}
For any non-zero $\mu\in {\cal M}_F(E)$,
\begin{equation*}
\lim_{r\to 0}r^{-\epsilon_0}\P_{\mu}(0<W\le r)
=e^{-\langle v,\,\mu\rangle}A(\psi(1))\langle v\phi_0^*,\,\mu\rangle/\Gamma(\epsilon_0+1),
\end{equation*}
where $\Gamma(\cdot)$ is the usual $\Gamma$ function, $\phi_0^*$ is
an egenfunction of $T^*_t$ corresponding to the eigenvalue $e^{\lambda^*_0t}$,
the operator $A$ is defined in \eqref{def-A}.
Furthermore
\begin{equation*}
\lim_{r\to\infty}r\widetilde{L}(r)^{-1}\P_{\mu}(W> r)=0.\end{equation*}
\end{thrm}

\begin{remark}{\rm
For a Galton-Watson process,
the small value probability problem of $W$ can be divided into two cases:
the Schr\"oder case and B\"ottcher case,  see \cite{Jones, MO}.
Suppose $\{Z_n, n\ge 0\}$ is a Galton-Watson process with offspring distribution
$\{p_n,n\ge 0\}$. Let $q$ be its extinction probability, $f(s)$ be the probability generating
function of $\{p_n,n\ge 0\}$, and $m>1$ be the mean of the offspring distribution.
Let $\gamma=f'(q)$.
\begin{itemize}
  \item [(1)] If $p_0+p_1>0$, then $F(s):=\lim_{n\to\infty}\gamma^{-n}(f_n(s)-q)$ exists,
  and  $F$ satisfies the Schr\"oder equation: $F(f(s))=\gamma f(s)$.
  Let $\epsilon=-\log \gamma/\log m$, then
  \begin{equation*}
  P(W\le r)\asymp r^{-\epsilon}.
  \end{equation*}
  \item [(2)] If $p_0+p_1=0$, then $\lambda=\min\{n:p_n>0\}\ge2$. In this case
  $G(s):=\lim_{n\to\infty}-\lambda^{-n}\log f_n(s)$ exists, and the function
  $\overline{G}=e^{-G}$ satisfies the B\"ottcher equation $\overline{G}(f)=\overline{G}^\lambda$. Let $\beta=\log \lambda/\log m$, then one can obtain
  \begin{equation*}
  -\log P(W\le r)\asymp r^{-\beta/(1-\beta)}.
  \end{equation*}
\end{itemize}

For the branching Markov process in \cite{Hering} and the superprocess in this paper,
the small value probability problem of the strong limit $W$ has only one case, the Schr\"oder case.
In fact, for the  branching Markov process $\{Z_t,t\ge0\}$ in \cite{Hering},
when the extinction probability is $0$, one can show that $\lim_{t\to\infty}e^{-\lambda^*_0t}\overline{F}_tf$ exists,
where $\overline{F}_tf:=P_{\delta_x}(e^{\langle\log f,Z_t\rangle})$.

Suppose that there exists $N>0$ such that the offspring distribution $\{p_n\}$ of
the Galton-Watson process satisfies $p_{n}=0$ for all $n\ge N$, then
\cite{BB} obtained the rate at which the tail probability of  $W$ tends to $0$.
For results on the rate at which tail probability of $W$ tends to $0$ for  multitype  Galton-Watson processes, see \cite{Jones}.
For superprocesses, under some condition, the rate at which tail probability of $W$
tends to $0$ is determined by the skeleton process (a branching Markov process)
of  $X$. When the branching mechanisms $n(x,{\rm d}r)$ is not $0$,
the offspring distribution $\{p_n\}$ (see \eqref{def-p2} and \eqref{def-pn}) of the skeleton
process of  $X$ does not satisfy this condition, thus we could not get the rate at which the tail probability $\P_{\mu}(W> r)$ tends to $0$ as $r\to\infty$.  We only obtain a weaker result. }
\end{remark}

In Section \ref{sec:bmp},
we will prove that $W$ is a compound Poisson random variable of the form $W=\sum_{n=1}^{N} Y_n$, where $N$ is a Poisson random variable,
$Y_j, j\ge 1$ is a sequence of independent and identically distributed random variables independent of  $N$.
We will also prove that the distribution of $Y_1$
is the distribution of the corresponding strong limit
of some branching Markov process.
In Section \ref{sec:proof},
we will analyze and estimate the Laplace transform and
characteristic function of $Y_1$, and show that
$Y_1$ restricted to $(0, \infty)$ has a density function,
thus proving \ref{main-dens}.
By using the Tauberian theorem, we can prove Theorem \ref{main-tail}.

\section{Compound Poisson random variable and branching Markov process}\label{sec:bmp}
\subsection{Compound Poisson random variable}

The Laplace exponent of $W$ is defined as
\begin{equation}\label{def-Phi}
\Phi(\theta,x):=-\log \P_{\delta_x}\exp\{-\theta W\}.
\end{equation}
Using the Markov property and the branching property, we have shown in \cite[(5.3)]{RSZ6} that
\begin{equation}\label{Phi}
\Phi(\theta,x)=V_t(\Phi(\theta e^{-\lambda_0t}, \cdot))(x).
\end{equation}

\begin{lemma}\label{infi-divis}
For any $x\in E$, there exists a finite measure $\pi(x,dr)$ on
$(0,\infty)$ such that $\pi(x,(0,\infty))=v(x)$, and
\begin{equation*}
\Phi(\theta,x)=\int_{(0,\infty)}\Big(1-e^{-\theta r}\Big)\pi(x,dr).\end{equation*}
\end{lemma}

\noindent\textbf{Proof:}
Since
\begin{equation*}
\P_{\mu}\big[e^{-\theta W}\big]=\Big(\P_{\mu/n}\big[e^{-\theta W}\big]\Big)^n,
\end{equation*}
the distribution of the random variable $W$ under $\P_{\mu}$ is infinitely divisible.
Since $W$ is non-negative, there exist a non-negative function $a(x)$ and
a $\sigma$-finite measure $\pi(x,dr)$ satisfying the condition
\begin{equation*}
\int_0^\infty (1\wedge r)\pi(x,dr)<\infty
\end{equation*}
such that
\begin{equation}\label{levy-ki}
  \Phi(\theta,x)=a(x)\theta+\int_{(0,\infty)}\Big(1-e^{-\theta r}\Big)\pi(x,dr).
\end{equation}
It follows from \cite[Theorem 1.2]{RSZ6} that
\begin{equation*}\Phi(\infty,x)=-\log \P_{\delta_x}(W=0)=v(x).\end{equation*}
From this one gets $a(x)=0$ and
\begin{equation*}\pi(x,(0,\infty))=v(x).\end{equation*}
The assertion of the lemma follows immediately.
\hfill $\Box$

\bigskip

It follows from Lemma \ref{infi-divis} that, under $\P_{\mu}$, $W$ is a compound Poisson
random variable, that is,
\begin{equation*}W=\sum_{n=1}^{N} Y_n,\end{equation*}
where $N$ is a Poisson random variable with parameter $\langle v,\mu\rangle$, $\{Y_j,j\ge1\}$
is a sequence of independent and identically distributed random variables with common
distribution $\frac{\int_E \pi(x,{\rm d}y)\, \mu(dx)}{\langle v,\mu\rangle}$ and independent of
$N$. From now on, we assume that $Y$ is a random variable with
distribution $\frac{\int_E \pi(x,{\rm d}y)\, \mu(dx)}{\langle v,\mu\rangle}$.

\begin{lemma}\label{density2}
For any non-zero $\mu\in {\cal M}_F(E)$, under $\P_{\mu}$ÏÂ,
the random variable $W$ restricted to $(0,\infty)$ has a density if and only if the random variable $Y$ has a density. Furthermore, if the density of $Y$ is $g_{\mu}(y)$, then for any  $0<a<b$,
\begin{equation*}\P_{\mu}(W\in(a,b))=\int _a^b f_{\mu}(y)\,d y,\end{equation*}
where
\begin{equation}\label{density}
  f_{\mu}(y)=\sum_{k=1}^\infty
  g_{\mu}^{*k}(y)
  \frac{\langle v,\mu\rangle^k}{k!}e^{-\langle v,\mu\rangle}.
\end{equation}
\end{lemma}

\begin{remark}\label{rem-den}
If for every $x\in E$, $Y$ has a density function $g(x,y)$ under
$\P_{\delta_x}$,  then
\begin{equation*}\pi(x, d y)=v(x)g(x,y)d y.\end{equation*}
Thus for every
$\mu\in {\cal M}_F(E)$, $Y$ has a density function under $\P_{\mu}$:
\begin{equation*}g_{\mu}(y)=\frac{\int_E v(x)g(x,y) \, \mu(d x)}{\langle v,\mu\rangle},\quad y>0.\end{equation*}
\end{remark}

It follows from Lemma \ref{density2}, it suffices to prove that the random variable $Y$
has a density function. For this we will analyze the Laplace transform and characteristic function
of $Y$. Define
\begin{equation*}
\psi(\theta,x):=\frac{v(x)-\Phi(\theta,x)}{v(x)}=v(x)^{-1}\int_{(0,\infty)}e^{-\theta r}\pi(x,dr),\quad \theta\ge 0.
\end{equation*}
Thus
\begin{equation*}
\P_{\mu} \Big(e^{-\theta Y}\Big)=
\frac{\langle v \psi(\theta,\cdot),\mu\rangle}{\langle v,\mu\rangle},\quad \theta\ge0.
\end{equation*}
Note that $\psi(\theta,x)$ is the Laplace transform of the distribution
$v(x)^{-1}\pi(x, dr)$. For any $x\in E$, $\theta\ge 0$, $\psi(\theta,x)\in [0, 1]$.

Similarly, we define
\begin{equation*}
\psi(i\theta,x):=\P_{\delta_x}\Big(e^{-i\theta Y}\Big)=v(x)^{-1}\int_{(0,\infty)}e^{-i\theta r}\pi(x,dr),\quad \theta\in \R.
\end{equation*}

For any $a>0$, let $\mathcal{D}_a:=\{f\in \mathcal{B}(E): 0\le f(x)\le a, x\in E\}$.
Define an operator $\overline{V}_t:\mathcal{D}_1\to \mathcal{D}_1$ by
\begin{equation}\label{barV}
 \overline{V}_tf(x):=
  \frac{v(x)-V_t(v(1-f))(x)}{v(x)},\qquad f\in
  \mathcal{D}_1.
\end{equation}
It follows from \eqref{Phi} that, for all $\theta\ge 0$,
\begin{equation}\label{psi}
  \psi(\theta,x)=\overline{V}_t(\psi(\theta e^{-\lambda_0t}, \cdot))(x).
\end{equation}
Obviously, we can extend the definition of $\overline{V}_t$ to the space of
all complex-valued functions on $E$
with sup norm less than or equal to 1.

In the next subsection, we will show that $\overline{V}_t$ is the Laplace functional of some
branching Markov process.   A skeleton decomposition of superprocesses was established under some conditions in \cite{CRY,EKW}. When the conditions of \cite{CRY,EKW} are satisfied, the
branching Markov process we are going to introduce below is just the skeleton process of the
$(\xi,\varphi)$-superprocess. \cite{EKW} dealt with the skeleton decomposition of super-diffusions,
while \cite{CRY} dealt with the skeleton decomposition of superprocesses with a symmetric
spatial motion. \cite{CRY,EKW} can not completely cover the superprocesses dealt with in this paper. In this paper we do not use the skeleton decomposition of superprocesses. We start with
that $W$ is a compound Poisson random variable, and introduce the corresponding
branching Markov process.

\subsection{Branching Markov processes}

Define
\begin{equation*}
N_t:=\frac{v(\xi_t)}{v(\xi_0)}\exp\Big\{-\int_0^t \frac{\varphi(\xi_s,v(\xi_s))}{v(\xi_s)}\,ds\Big\}.
\end{equation*}

\begin{lemma}
Under $\Pi_x$,
$\{N_t:t\ge0\}$ is a non-negative martingale with respect to the filtration
$\{\mathcal{H}_t, t\ge 0\}$,
and $\Pi_x(N_t)=1$.
\end{lemma}

\noindent\textbf{Proof:}
It follows from the Markov property and the branching property that, for any $t>0$, $v(x)=V_tv(x)$.
Thus
\begin{equation}\label{eq:v}
  v(x)+\Pi_x\int_0^t\varphi(\xi_s, v(\xi_s)){\rm d}s=\Pi_x v(\xi_t).
\end{equation}
It is easy to see that
\begin{equation*}|\varphi(x,z)|\le 2M(z+z^2),\quad z\ge0,\end{equation*}
Thus
\begin{equation}\label{bound-varphi}
  \frac{|\varphi(x,v(x))|}{v(x)}\le 2M(1+v(x))\le 2M(1+\|v\|_\infty).
\end{equation}
Hence it follows from the Feynman-Kac formula that
\begin{equation}\label{eq:v2}
  v(x)=\Pi_x\Big[ \exp\Big\{-\int_0^t \frac{\varphi(\xi_s,v(\xi_s))}{v(\xi_s)}\,ds\Big\}v(\xi_t)\Big].
\end{equation}
It follows immediately from the Markov property and \eqref{eq:v2} that
\begin{align*}
 \Pi_x(N_{t+s}|\mathcal{H}_t)
 =&v(\xi_0)^{-1}\exp\Big\{-\int_0^t \frac{\varphi(\xi_s,v(\xi_s))}{v(\xi_s)}\,ds\Big\}\\
 &\times\Pi_{\xi_t}\Big[v(\xi_s)\exp\Big\{-\int_0^s \frac{\varphi(\xi_s,v(\xi_s))}{v(\xi_s)}\,ds\Big\}
 =N_t.
\end{align*}
Thus $\{N_t:t\ge0\}$ is a non-negative martingale.
\hfill$\Box$

We use the martingale $\{N_t\}$ to define a new probability measure $\overline{\Pi}_x$:
\begin{equation*}
\frac{d\overline{\Pi}_x}{d\Pi_x}\Big|_{\mathcal{H}_t}=N_t,\qquad t\ge 0.
\end{equation*}

\begin{prop}\label{prop:barV}
For $f\in \mathcal{D}_1$,
\begin{equation}\label{bar-V2}
   \overline{V}_tf(x)=\overline{\Pi}_x\int_0^t\varphi^*(\xi_s,\overline{V}_{t-s}f(\xi_s)){\rm d}s+\overline{\Pi}_x f(\xi_t),
\end{equation}
where
\begin{equation*}\varphi^*(x,\lambda):=\frac{\varphi(x,v(x)(1-\lambda))-\varphi(x,v(x))(1-\lambda)}{v(x)},\quad \lambda\in[0,1].\end{equation*}
\end{prop}

\noindent\textbf{Proof:}
It follows from \eqref{V_t} that
for $f\in \mathcal{D}_1$,
\begin{equation*}V_t(v(1-f))(x)+\Pi_x\int_0^t\varphi(\xi_s, V_{t-s}(v(1-f))(\xi_s)){\rm d}s=\Pi_x v(\xi_t)(1-f(\xi_t)).\end{equation*}
Thus by \eqref{eq:v}, we have
\begin{equation*}
v(x)\overline{V}_tf(x)+\Pi_x\int_0^t\varphi(\xi_s, v(\xi_s))-\varphi(\xi_s, v(\xi_s)(1-\overline{V}_{t-s}f(\xi_s))){\rm d}s
=\Pi_x v(\xi_t)f(\xi_t).
\end{equation*}
Hence
\begin{align*}
  v(x)\overline{V}_tf(x)=&\Pi_x\int_0^tv(\xi_s)\varphi^*(\xi_s,\overline{V}_{t-s}f(\xi_s)){\rm d}s\\
  &-\Pi_x\int_0^t\frac{\varphi(\xi_s,v(\xi_s))}{v(\xi_s)}v(\xi_s)\overline{V}_{t-s}f(\xi_s){\rm d}s+\Pi_x v(\xi_t)f(\xi_t).
\end{align*}
It follows from the Feynman-Kac  formula that
\begin{equation*}
   v(x)\overline{V}_tf(x)
   =\Pi_x\int_0^te^{-\int_0^s\frac{\varphi(\xi_u,v(\xi_u))}{v(\xi_u)}\,du}v(\xi_s)\varphi^*(\xi_s,\overline{V}_{t-s}f(\xi_s)){\rm d}s
   +{\Pi}_x\Big[e^{-\int_0^t\frac{\varphi(\xi_s,v(\xi_s))}{v(\xi_s)}\,ds}v(\xi_t) f(\xi_t)\Big],
\end{equation*}
from which \eqref{bar-V2} follows immediately
\hfill$\Box$

With the preparation above, we now introduce a branching Markov process corresponding to
$\overline{V}_t$. By the definition of $\varphi^*(x,\lambda)$, we have
\begin{align*}
  \varphi^*(x,\lambda)=&\frac{\varphi(x,v(x)(1-\lambda))-\varphi(x,v(x))(1-\lambda)}{v(x)} \\
  =&\beta(x)v(x)(\lambda^2-\lambda)+v(x)^{-1}\int_0^\infty\Big((e^{r\lambda v(x)}-1+\lambda)e^{-rv(x)}-\lambda\Big)n(x,dr)\\
  =&\beta(x)v(x)\lambda^2+\sum_{n=2}^\infty\int_0^\infty\frac{v(x)^{n-1}(r\lambda )^n}{n!}e^{-rv(x)}\,n(x,dr)\\
  &-\lambda\Big(\beta(x)v(x)+v(x)^{-1}\int_0^\infty(e^{rv(x)}-1-rv(x))e^{-rv(x)}\,n(x,dr)\Big).
\end{align*}
Thus we have
\begin{equation}\label{2.1}
  \varphi^*(x,\lambda)=
  b(x)\Big(\sum_{n=2}^\infty \lambda^np_n(x)-\lambda\Big),
\end{equation}
where
\begin{equation*}
b(x)=\beta(x)v(x)+v(x)^{-1}\int_0^\infty(e^{rv(x)}-1-rv(x))e^{-rv(x)}\,n(x,dr);
\end{equation*}
\begin{align}
p_2(x) =& \frac{v(x)}{b(x)}
  \Big(\beta(x)+\frac{1}{2}\int_0^\infty r^2e^{-v(x)r}n(x,dr)\Big); \label{def-p2}\\
  p_n(x) =& \frac{v^{n-1}(x)}{n!b(x)}
  \int_0^\infty r^ne^{-v(x)r}n(x,dr),\quad n>2.\label{def-pn}
\end{align}
It is easy to verify that $\sum_{n=2}^\infty p_n(x)=1$ and $b(x)$ is a bounded non-negative function. In fact,
\begin{align*}
  b(x)
  \le& \beta(x)v(x)+v(x)^{-1}\int_0^\infty ((rv(x))\wedge (rv(x))^2)\,n(x,dr)\\
  \le& \beta(x)v(x) +v(x)\int_0^1 r^2\,n(x,dr)+\int_1^\infty r\,n(x,dr)\le M\|v\|_\infty+M.
\end{align*}
It is also easy to see that $b(x)>0$.

Consider a branching Markov process $\{Z_t, t\ge 0; P_{\nu}\}$ with spatial motion
$\{\xi_t, t\ge 0; \overline{\Pi}_x\}$, branching rate function $b(x)$ and spatially dependent offspring distribution $\{p_n(x):n\ge2\}$. Then for any $g\in\mathcal{B}_b^+(E)$,
\begin{equation*}P_{\delta_x}(e^{-\langle g,Z_t\rangle})=\overline{V}_t(e^{-g})(x),\end{equation*}
and
\begin{equation}\label{smg-bp}
  Q_tg(x):=P_{\delta_x}(\langle g,Z_t\rangle)
  =\overline{\Pi}_x\left(\exp\Big\{\int_0^t\frac{\partial}{\partial \lambda}\varphi^*(\xi_s,1)\,ds\Big\}g(\xi_t)\right)=v(x)^{-1}T_t(vg)(x),
\end{equation}
where the last equality follows from the definitions of $\varphi^*$ and $\overline{\Pi}_x$.
Hence the first eigenvalue of the infinitesimal generator of the semigroup $\{Q_t\}$ is
$\lambda_0$,
and $v(x)^{-1}\phi_0(x)$ is the corresponding eigenfunction. It follows from the
boundedness of $v$ that $\{Q_t\}$ is also intrinsically untracontractive, and thus
condition (M) in \cite{Hering} holds. Hence it follows from \cite[Proposition 3.6]{Hering} that
there exist a non-negative function $\overline{\gamma}_t$ and a non-degenerate random variable
$W^Z$ such that
\begin{equation*}
\overline{\gamma}_t\langle v^{-1}\phi_0,Z_t\rangle\to W^Z,\quad P_{\nu}\mbox{-a.s.},
\end{equation*}
and the Laplace transform of $W^Z$, defined by
\begin{equation*}
\psi^Z(\theta,x):=P_{\delta_x}\Big(e^{-\theta W^Z}\Big),\quad \theta\in \R,
\end{equation*}
is a solution of \eqref{psi}. We already know that the Laplace transform $\psi(\theta,x)$ of $Y$
is also a solution of \eqref{psi}, thus it follows from \cite[Proposition 3.8]{Hering} that
there exists $a\in(0,\infty)$ such that $(Y,\P_{\delta_x})$ and $(aW^Z,P_{\delta_x})$
have the same distribution. Since $p_0(x)=p_1(x)=0$, the extinction probability of $Z$
is 0. Using  the assertions about $W^Z$ in \cite[Propositions 5.1, 5.10 and 5.11]{Hering},
one can deduce the corresponding properties of $(Y,\P_{\delta_x})$, and thus obtaining the
proofs of Theorem \ref{main-dens} and Theorem \ref{main-tail}.

In Theorem \ref{main-tail}, the semigroup $\{T_t^*\}$
(see Subsection \ref{sec:est}),
especially the first eigenvalue $\lambda_0^*$ of its infinitesimal generator
and its corresponding eigenfunction,  play very important roles. Theorem \ref{main-tail}
contains another important operator $A$, which is determined by the limit of
$e^{-\lambda_0^*t}\overline{V}_tf(x)$, see \eqref{def-A}. The semigroup
$\{\delta\bar{F}_t(0),t\ge0\}$ in \cite{Hering} coincides with the semigroup
$\{T_t^*\}$ of this paper, and the operator $Q$ there coincides with our operator $A$,
but \cite{Hering} did not give explicit expressions for these two quantities.
For completeness, we do not quote the conclusions of \cite{Hering} directly.
In Subsection \ref{sec:est}, we will give the definitions of $\{T_t^*\}$ and $A$.
In Subsections \ref{3.2} and \ref{3.3}, we will give the proofs of Theorem \ref{main-dens} and Theorem \ref{main-tail}. The main ideas are similar to that of \cite{Hering}.

\section{Proofs of Main Results}\label{sec:proof}

\subsection{Estimates on the operator $\overline{V}_t$}\label{sec:est}

In this subsection,
we will give some estimates on the operator $\overline{V}_t$.
We then use these estimates and \eqref{psi} to obtain some estimates on the Laplace transform $\psi(\theta,x)$. In the proof below, $C$ stands for a constant whose value
might change from one appearance to another.

We first list some estimates from \cite{RSZ6} that we will use in this paper.

\begin{itemize}
  \item [(1)] {\bf{Estimates on the semigroup $\{T_t\}$}:}
  It follows from \cite[Theorem 2.7]{KiSo08b} that,  under Assumptions \ref{assum1'}--\ref{assum63},
  for any $\delta>0$, there exist constants $\gamma=\gamma(\delta)>0$ and $c=c(\delta)>0$ such that for all
  $(t,x,y)\in [\delta,\infty)\times E\times E$, we have
  \begin{equation}\label{density0}
   \left|e^{-\lambda_0t}q(t,x,y)-\phi_0(x)\psi_0(y)\right|\le ce^{-\gamma t}\phi_0(x)\psi_0(y).
  \end{equation}
   Take $t$ large enough so that $ce^{-\gamma t}<\frac12$, then
  \begin{equation*}
  e^{-\lambda_0 t}q(t, x, y)\ge \frac12\phi_0(x)\psi_0(y).
  \end{equation*}
  Since $q(t, x, \cdot)\in L^1(E, m)$, we have $\psi_0\in L^1(E, m)$.
  Thus for any $f\in\mathcal{B}_b^+(E)$,  we have $\la f,\psi_0\ra_m<\infty$.
  Consequently, for any $f\in\mathcal{B}_b^+(E)$¼°$(t,x)\in[\delta,\infty)\times E$, we have
  \begin{equation}\label{IU-T_t}
  \left|e^{-\lambda_0t}T_tf(x)-\langle f,\psi_0\rangle_m \phi_0(x)\right|
  \le ce^{-\gamma t}\langle |f|,\psi_0\rangle_m \phi_0(x)
  \end{equation}
  and
  \begin{equation}\label{IU-T_t'}
   (1-ce^{-\gamma t})\langle |f|,\psi_0\rangle_m\phi_0(x)\le e^{-\lambda_0t}T_t|f|(x)\le  (1+c)\langle |f|,\psi_0\rangle_m\phi_0(x).
  \end{equation}

    \item [(2)]{\bf{$v$ and $\phi_0$ are comparable:}}
    It follows from \cite[Lemma 4.4]{RSZ6} that
  \begin{equation}\label{dom-v-below}
  v(x)=V_1(v)(x)\ge C T_1(v)(x)\ge C\phi_0(x).
  \end{equation}
  Furthermore,
  \begin{equation}\label{dom-v-above}v(x)=V_1v(x)\le T_1v(x)\le C \phi_0(x).\end{equation}

     \item [(3)]
       It follows from \cite[Proposition 5.3 and Lemma 4.6]{RSZ6} that
  \begin{equation*}
  \gamma_t=\langle \Phi(e^{-\lambda_0t}),\psi_0 \rangle_m,\end{equation*}
  \begin{equation}\label{limit-Phi}
  \Phi(e^{-\lambda_0t},x)=(1+h_t(x))e^{-\lambda_0t}L(t)\phi_0(x),
  \end{equation}
  where
  $\lim_{t\to\infty}\|h_t\|_\infty=0.$
  \end{itemize}

Define a  semigroup ${P}^{\varphi'}_t$ by
\begin{equation}\label{T_t*}
  {P}^{\varphi'}_tf(x)
  :=\Pi_x\Big(f(\xi_t)e^{-\int_0^t \partial_\lambda\varphi(\xi_s,v(\xi_s))\,ds}\Big).
\end{equation}

It follows from \eqref{e:M} and the boundedness of $v$ that $\partial_\lambda\varphi(x,v(x))$
is bounded. By using the same argument as in the paragraph above \eqref{6condition2}, one can show that the semigroup $({P}^{\varphi'}_t)$ is also intrinsically ultracontractive.
Let $\lambda^*_0$ be the largest (simple) eigenvalue of the infinitesimal generator of $({P}^{\varphi'}_t)$, let  $\overline{\phi}_0$ and  $\overline{\psi}_0$ be, respectively, eigenfunctions  of the infinitesimal generators of $({P}^{\varphi'}_t)$ and its dual semigroup
corresponding to $\lambda^*_0$. $\overline{\phi}_0$ and  $\overline{\psi}_0$ can be chosen
to be strictly positive continuous functions on $E$ and satisfy $\|\overline{\phi}_0\|_2=1$,
$\langle \overline{\phi}_0,\overline{\psi}_0\rangle_m=1$.  Furthermore, $\overline{\phi}_0$
is a bounded function, and $\overline{\psi}_0\in L^1(E,m)$.

 It follows from \cite[Theorem 2.7]{KiSo08b} that,
 under Assumptions \ref{assum1'}--\ref{assum63},
  for any $\delta>0$, there exist constants $\gamma=\gamma(\delta)>0$ and $c=c(\delta)>0$ such that for all
  $(t,x,y)\in [\delta,\infty)\times E\times E$, we have
  \begin{equation}\label{T_t^*3}
  \left|e^{-\lambda_0^*t}P_t^{\varphi'}f(x)-\langle f,\overline{\psi}_0\rangle_m \overline{\phi}_0(x)\right|
  \le ce^{-\gamma t}\langle |f|,\overline{\psi}_0\rangle_m \overline{\phi}_0(x).
\end{equation}

Based on this, we define another semigroup
\begin{equation*}
T_t^*f(x):=v(x)^{-1}{P}^{\varphi'}_t(vf)(x)
=v(x)^{-1}\Pi_x\Big((vf)(\xi_t)e^{-\int_0^t \partial_\lambda\varphi(\xi_s,v(\xi_s))\,ds}\Big)
=\overline{\Pi}_x(f(\xi_t)e^{-\int_0^t b(\xi_s)\,ds}).
\end{equation*}
Let $\phi^*_0(x):=v(x)^{-1}\overline{\phi}_0(x)$ and $\psi^*_0(x):=v(x)\overline{\psi}_0(x)$.
It follows from \eqref{T_t^*3} that
\begin{equation}\label{T*}
  \left|e^{-\lambda_0^*t}T_t^{*}f(x)-\langle f,{\psi}^*_0\rangle_m {\phi}_0^*(x)\right|
  \le ce^{-\gamma t}\langle |f|,{\psi}^*_0\rangle_m {\phi}^*_0(x).
\end{equation}
Note that $\partial_\lambda\varphi(x,\lambda)\ge -\alpha(x)$. Hence
\begin{equation*}
\overline{\phi}_0(x)=e^{\lambda_0^*}{P}^{\varphi'}_1(\overline{\phi}_0)(x)
\le CT_1(\overline{\phi}_0)(x)\le C\phi_0(x).
\end{equation*}
Using \eqref{dom-v-below} we see that $\|\phi^*_0\|_\infty<\infty$. It is also easy to see that
$\psi^*_0\in L^1(E,m)$.

It can be shown that the semigroup $\{T_t^*,t\ge 0\}$ defined above coincides with
the semigroup $\{\delta \bar{F}_t(0)\}$ defined in \cite{Hering}, where $\{\delta \bar{F}_t(0)\}$ is defined via a Fr\'{e}chet derivative.

\begin{lemma}
\begin{equation*}\lambda_0^*<0.\end{equation*}
\end{lemma}

\noindent\textbf{Proof:}
It follows from $V_tv(x)=v(x)$ and \cite[Lemma 4.1]{RSZ8} that
\begin{equation*}
\P_{\delta_x}\Big(\langle f,X_t\rangle e^{-\langle v,X_t\rangle}\Big)
=\Pi_x\Big(f(\xi_t)e^{-\int_0^t \partial_z\varphi(\xi_s,v(\xi_s))\,ds}\Big)e^{-v(x)}.
\end{equation*}
Thus,
\begin{equation}\label{T*2}
  T^{*}_tf(x)=v(x)^{-1}e^{v(x)}\P_{\delta_x}\Big(\langle vf,X_t\rangle e^{-\langle v,X_t\rangle}\Big).
\end{equation}
By \cite[Lemma 3.2]{RSZ6}, we have
\begin{equation*}
\P_{\delta_x}(\lim_{t\to\infty}\langle v,X_t\rangle=0)
=1-\P_{\delta_x}(\lim_{t\to\infty}\langle v,X_t\rangle=\infty)=e^{-v(x)}.
\end{equation*}
Hence by the dominated convergence theorem, we have
\begin{equation*}
\lim_{t\to\infty}T_t^*1(x)=0.
\end{equation*}
Combining with \eqref{T*}, we immediately get $\lambda_0^*<0.$
\hfill$\Box$

\begin{lemma}\label{barVto0}
For any  $a\in[0,1)$, there exists a constant $c(a)>0$ such that for $t\ge1$,
\begin{equation*}
\overline{V}_tf(x)\le c(a)e^{\lambda_0^*(1-a)t}{\phi}^*_0(x)^{1-a},
\quad \forall f\in {\cal D}_a.
\end{equation*}
\end{lemma}
\noindent\textbf{Proof:}
Note that
\begin{equation*}\sup_{0\le \lambda\le a}\frac{\varphi^*(x,\lambda)}{\lambda}=\sup_{0\le \lambda\le a}
b(x)(\sum_{n=2}^\infty p_n(x)\lambda^{n-1}-1)\le b(x)(a-1).\end{equation*}
Since $\varphi^*(x,\lambda)\le 0$ for all $\lambda\in [0, 1]$,
it follows from  \eqref{bar-V2} that for any
$f\in \mathcal{D}_1$,
\begin{equation}\label{11.5}
  \overline{V}_tf(x)\le \overline{\Pi}_xf(x)\le \|f\|_\infty.
\end{equation}
Thus for $f\in \mathcal{D}_a$,
\begin{equation}\label{11.4}
  \varphi^*(x,\overline{V}_{t-s}f(x))\le b(x)(a-1)\overline{V}_{t-s}f(x).
\end{equation}
It follows from \eqref{bar-V2} and the Feynman-Kac formula that, for $t\ge 1$,
\begin{align*}
 \displaystyle \overline{V}_tf(x)=&\displaystyle \overline{\Pi}_x\int_0^t  e^{-(1-a)\int_0^s
 b(\xi_u)\,du}\Big[\varphi^*(\xi_s,\overline{V}_{t-s}f(\xi_s))+(1-a)b(\xi_s)\overline{V}_{t-s}f(\xi_s)\Big]{\rm d}s\\
 &\displaystyle +\overline{\Pi}_x \Big[e^{-(1-a)\int_0^t b(\xi_s)\,ds}f(\xi_t)\Big]\\
 \le&\displaystyle  a\overline{\Pi}_x \Big[e^{-(1-a)\int_0^t b(\xi_s)\,ds}\Big]
 \le\displaystyle  a\left[\overline{\Pi}_x \Big[e^{-\int_0^t b(\xi_s)\,ds}\Big]\right]^{1-a}\\
 =&\displaystyle  a\left[T^*_t1(x)\right]^{1-a}
 \le a(1+c)e^{\lambda^*_0(1-a)t}\langle 1,{\psi}^*_0\rangle_m^{1-a} {\phi}^*_0(x)^{1-a},
\end{align*}
where the last inequality follows from \eqref{T*}.
\hfill$\Box$

\begin{lemma}\label{lemV}
For any $f\in {\cal D}_1$,
\begin{equation*}
T_{t}^*f(x)\le \overline{V}_tf(x)\le (1+\|f\|_\infty
e^{\|b\|_\infty t})T_{t}^*f(x).
\end{equation*}
\end{lemma}
\noindent\textbf{Proof:}
It follows from \eqref{11.5} that
\begin{equation}\label{11.6}
  \overline{V}_tf(x)\le \overline{\Pi}_xf(x)\le
  e^{\|b\|_\infty t}T_t^*f(x).
\end{equation}
Using \eqref{bar-V2} and the Feynman-Kac formula, we can get
\begin{equation}\label{bar-V3}
   \overline{V}_tf(x)=\int_0^tT_s^*[\varphi^*_0(\cdot,\overline{V}_{t-s}f)](x){\rm d}s+T^*_t (f)(x),
\end{equation}
where
$\varphi_0^*(x,\lambda)=\varphi^*(x,\lambda)+
b(x)\lambda\ge 0.$
Hence
 $\overline{V}_tf(x)\ge T_{t}^*f(x).$
Note that
\begin{equation}\label{11.11}
  \varphi_0^*(x,\lambda)\le
  b(x)\lambda^2, \qquad 0\le \lambda\le 1.
\end{equation}
Combining \eqref{11.5}  and \eqref{11.6}, we get
\begin{equation*}\varphi^*_0(x,\overline{V}_{t-s}f(x))\le
b(x)\overline{V}_{t-s}f(x)^2\le \|f\|_\infty\|b\|_\infty e^{\|b\|_\infty
(t-s)}T_{t-s}^*f(x).\end{equation*}
Hence,
\begin{equation*}
\int_0^tT_s^*[\varphi^*_0(\cdot,\overline{V}_{t-s}f)](x){\rm d}s\le  \|f\|_\infty\int_0^t
\|b\|_\infty e^{\|b\|_\infty (t-s)}{\rm d}sT_{t}^*f(x)
\le \|f\|_\infty e^{\|b\|_\infty t}T_{t}^*f(x).
\end{equation*}
Summarizing the above, we get
\begin{equation*} T_{t}^*f(x)\le \overline{V}_tf(x)\le (1+\|f\|_\infty
e^{\|b\|_\infty t})T_{t}^*f(x).
\end{equation*}
\hfill$\Box$

For any $f\in \mathcal{D}_1$,
define
\begin{equation}\label{def-A}
A(f):=
\int_{0}^\infty e^{-\lambda_0^* s}\langle \varphi^*_0(\cdot,\overline{V}_{s}f),\psi_0^*\rangle_m{\rm d}s
+\langle f,\psi_0^*\rangle_m.
\end{equation}

\begin{lemma}\label{barV4}
For any $a\in[0,1)$ and
$f\in \mathcal{D}_a$,
\begin{equation*}
\sup_{t>0}e^{-\lambda_0^* t}\|\overline{V}_tf\|_\infty<\infty.
\end{equation*}
Furthermore,
\begin{equation*}
\lim_{t\to\infty} e^{-\lambda_0^* t}\overline{V}_tf(x)
=A(f)\phi_0^*(x),
\end{equation*}
where $A(f)$is defined in \eqref{def-A}.
\end{lemma}
\noindent\textbf{Proof:}
Note that $\phi^*_0(x)$ is bounded.
It follows from Lemma \ref{barVto0} that there exists
$s_0>1$ such that
\begin{equation*}
\overline{V}_{s_0}f(x)\le 1/4, \quad \forall f\in {\cal D}_a.
\end{equation*}
Using this and Lemma \ref{barVto0} we  obtain that, for any $s>s_0+1$,
\begin{equation}\label{11.10}
\overline{V}_sf(x)=\overline{V}_{s-s_0}(\overline{V}_{s_0}f)(x)
\le\overline{V}_{s-s_0}(1/4)(x)
\le Ce^{3\lambda^*_0s/4} \phi^*_0(x)^{3/4}.
\end{equation}
It follows from \eqref{bar-V3} that for any $t>s_0+1$,
\begin{align}\label{bar-V4}
   &e^{-\lambda_0^* t}\overline{V}_tf(x)=\int_0^te^{-\lambda_0^* s}e^{-\lambda_0^* (t-s)}T_{t-s}^*[\varphi^*_0(\cdot,\overline{V}_{s}f)](x){\rm d}s+e^{-\lambda_0^* t}T^*_t (f)(x)\nonumber\\
   =&\Big(\int_0^{s_0+1}+\int_{s_0+1}^t\Big)e^{-\lambda_0^* s}e^{-\lambda_0^* (t-s)}T_{t-s}^*[\varphi^*_0(\cdot,\overline{V}_{s}f)](x){\rm d}s+e^{-\lambda_0^* t}T^*_t (f)(x)\nonumber\\
   =:&J_1(t,x)+J_2(t,x)+J_3(t,x).
\end{align}
For $J_3$, using  \eqref{T*} we can easily get
$J_3(t,x)\le C\langle f,\psi_0^*\rangle_m\phi_0^*(x)$ and
\begin{equation*}\lim_{t\to\infty}J_3(t,x)=\langle f,\psi_0^*\rangle_m\phi_0^*(x).\end{equation*}

For $J_2(t,x)$, we can use \eqref{11.10} and \eqref{11.11} to get that,  for any $t>s>s_0+1$,
\begin{equation*}
e^{-\lambda_0^* (t-s)}|T_{t-s}^*[\varphi^*_0(\cdot,\overline{V}_{s}f)](x)|\le Ce^{3\lambda^*_0s/2}
e^{-\lambda_0^* (t-s)}T_{t-s}^*[(\phi^*_0)^{3/2}](x)\le Ce^{3\lambda^*_0s/2}\phi^*_0(x).
\end{equation*}
Hence,
\begin{equation*}
|J_2(t,x)|\le C \int_{s_0+1}^te^{\lambda_0^* s/2}{\rm d}s\phi^*_0(x)\le C\phi^*_0(x).
\end{equation*}
It follows from the dominated convergence theorem that
\begin{equation*}
\lim_{t\to\infty}J_2(t,x)=\int_{s_0+1}^\infty e^{-\lambda_0^* s}\langle \varphi^*_0(\cdot,\overline{V}_{s}f),\psi_0^*\rangle_m{\rm d}s\phi^*_0(x).
\end{equation*}

Finally, we deal with $J_1(t,x)$. Since $\overline{V}_{s}f(x)\le 1$, we have $\varphi^*_0(x,\overline{V}_{s}f(x))\le \|b\|_\infty$.
Thus for $t-s>t-s_0>1$, we have
\begin{equation*}
e^{-\lambda_0^* (t-s)}|T_{t-s}^*[\varphi^*_0(\cdot,\overline{V}_{s}f)](x)|\le C e^{-\lambda_0^* (t-s)}T_{t-s}^*1(x)\le C\phi^*_0(x).
\end{equation*}
Hence
\begin{equation*}J_1(t,x)<C\phi^*_0(x),\end{equation*}
and it follows from the dominated convergence theorem that
\begin{equation*}
\lim_{t\to\infty}J_1(t,x)=\int_0^{s_0+1}e^{-\lambda_0^* s}\langle \varphi^*_0(\cdot,\overline{V}_{s}f),\psi_0^*\rangle_m{\rm d}s\phi^*_0(x).
\end{equation*}
Summarizing, we get the conclusion of the lemma.
\hfill$\Box$

\subsection{Proof of Theorem \ref{main-dens}}\label{3.2}

It follows from Lemma \ref{density2} that  to prove Theorem \ref{main-dens},
it suffices to show that $Y$ has a density $g_{\mu}(y)$ and that,
for any $y>0$, $g_{\mu}(y)>0$.
In this subsection,
we will show that  $Y$ has a density function by analyzing the
properties of the characteristic function of  $Y$.
By \eqref{psi}, we have
\begin{equation}\label{psi2}
  \psi(i\theta,x)=\overline{V}_t(\psi(i\theta e^{-\lambda_0t}, \cdot))(x),\quad \theta\in \R.
\end{equation}

For simplicity, for any $\theta\in \R$, we write $\psi(i\theta,\cdot)$
as $\psi(i\theta)$;  similarly, for any $\theta>0$, we write $\psi(\theta,\cdot)$
as $\psi(\theta)$.

\begin{lemma}\label{lem:11}
For any bounded closed interval $I$ not containing $0$, we have
\begin{equation*}
\sup_{\theta\in I}\|\psi(i\theta)\|_\infty<1.
\end{equation*}
\end{lemma}

\noindent \textbf{Proof:}
It is easy to see that
\begin{align*}
&|\|\psi(i\theta)\|_\infty-\|\psi(i(\theta+\epsilon))\|_\infty|\le \|\psi(i\theta)-\psi(i(\theta+\epsilon))\|_\infty\\
&=\|\overline{V}_1(\psi(i\theta e^{-\lambda_0}))-\overline{V}_1(\psi(i(\theta+\epsilon)e^{-\lambda_0}))\|_\infty.
\end{align*}
It is well known that, for any $|x_j|\le 1,$ $|y_j|\le 1$, it holds that
\begin{equation*}|\prod_{j=1}^n x_j-\prod_{j=1}^n y_j|\le \sum_{j=1}^n |x_j-y_j|.\end{equation*}
For any complex-valued function $f$ on $E$ with sup norm less or equal to 1, we have
\begin{equation*}
\overline{V}_tf(x)=P_{\delta_x}\prod_{u\in \mathcal{L}_t}f(\xi_t(u)),
\end{equation*}
where $\mathcal{L}_t$ is the collection of particles of the branching Markov process
$Z$ which are alive at time $t$, $\xi_t(u)$ stands for the position of particle $u$ at time $t$.
Thus,
\begin{align*}
&|\overline{V}_1(\psi(i\theta e^{-\lambda_0}))(x)-\overline{V}_1(\psi(i(\theta+\epsilon)e^{-\lambda_0}))(x)|\\
\le& P_{\delta_x}\langle|\psi(i\theta e^{-\lambda_0})-\psi(i(\theta+\epsilon)e^{-\lambda_0})|,Z_1\rangle\\
=&v(x)^{-1}T_1(v|\psi(i\theta e^{-\lambda_0}))-\psi(i(\theta+\epsilon)e^{-\lambda_0}))|)(x)\\
\le & C\langle |\psi(i\theta e^{-\lambda_0}))-\psi(i(\theta+\epsilon)e^{-\lambda_0}))|,\psi_0\rangle_m\to 0,\quad \epsilon \to0.
\end{align*}
The equality above is due to \eqref{smg-bp},
the last inequality is due to \eqref{11.6} and \eqref{dom-v-below}, and the last limit is due to the continuity of the characteristic function and the dominated convergence theorem.
Thus  $\|\psi(i\theta)\|_\infty$ is continuous in $\theta$. Now, we only need to show that,
for any $\theta\neq 0$, $\|\psi(i\theta)\|_\infty<1$.

We use contradiction. Suppose that for all $\theta\in \mathbb{R}$ and  $x\in E$,
$|\psi(i\theta,x)|=1$. Then by the uniqueness of characteristic functions,  there exists
a positive-valued function $c(x)$ such that $P_{\delta_x}(Y=c(x))=1$,
that is,  $\psi(\theta,x)=e^{-\theta c(x)}$.
Using \eqref{psi} (with $\theta$ replaced by $\theta e^{\lambda_0t}$),
\begin{equation*}\exp\{-\theta e^{\lambda_0t}c(x)\}=P_{\delta_x}(e^{-\theta \langle c,Z_t\rangle}),\end{equation*}
that is, the Laplace transform of the random variable $\langle c,Z_t\rangle$
is the same as $e^{\lambda_0t}c(x)$,
thus $P_{\delta_x}(\langle c,Z_t\rangle=e^{\lambda_0t}c(x))=1$.
By the definition of branching Markov processes, we know that $\langle c,Z_t\rangle$
can not be concentrated at one point, a contradiction!
Thus there exist $\theta_0\in\R$ and $x_0\in E$
such that $|\psi(i\theta_0,x_0)|<1$. Hence there exists $\delta=\delta(x_0)>0$, such that
$|\psi(i\theta,x_0)|<1$ for all $|\theta|\in (0,\delta)$.  Since $x\to \psi(i\theta,x)$ is a
continuous function, for all  $|\theta|\in (0,\delta)$, we have
\begin{equation*}m({y\in E:|\psi(i\theta,y)|<1})>0.\end{equation*}
For any complex-valued function $f$ on $E$ with sup norm less or equal to 1,
\begin{equation*}|\overline{V}_tf(x)|=|P_{\delta_x}\prod_{u\in \mathcal{L}_t}f(\xi_t(u))|\le P_{\delta_x}\prod_{u\in \mathcal{L}_t}|f|(\xi_t(u))=\overline{V}_t|f|(x).\end{equation*}
Thus by \eqref{11.6}, we have
\begin{equation*}
  1-|\psi(i\theta e^{\lambda_0t},x)|\ge 1-\overline{V}_t(|\psi(i\theta)|)(x) \ge \overline{\Pi}_x(1-|\psi(i\theta,\xi_t)|).
\end{equation*}
Note that
\begin{equation*}\frac{\varphi(x,v(x))}{v(x)}\le -\alpha(x)+\beta(x)v(x)+\frac{1}{2}v(x)\int_0^1 r^2n(x,dr)+\int_1^\infty r n(x,dr)
\le -\alpha(x)+M(\|v\|_\infty+1).\end{equation*}
Suppose that $c$ and $\gamma$ are the constants from \eqref{IU-T_t'}. For $t$ large enough,
we have $1-ce^{-\gamma t}>0$, hence by the definition of $\overline{\Pi}_x$,
\begin{align*}
&\overline{\Pi}_x(1-|\psi(i\theta,\xi_t)|)\ge v(x)^{-1}
e^{-M(\|v\|_\infty+1)t}T_t(v(1-|\psi(i\theta)|))(x)\\
\ge& e^{-M(\|v\|_\infty+1)t}(1-ce^{-\gamma t})
e^{\lambda_0t}\langle v(1-|\psi(i\theta)|),\psi_0\rangle_m\frac{\phi_0(x)}{v(x)},
\end{align*}
where the last inequality is due to \eqref{IU-T_t'}.

It follows from \eqref{dom-v-above} that,
for $|\theta|\in (0,\delta)$ and $t$ sufficiently large,
$\|\psi(i\theta e^{\lambda_0t})\|_\infty<1$. Thus  for all $\theta\neq 0$, $\|\psi(i\theta)\|_\infty<1.$

Summarizing, we get the conclusion of the lemma.
\hfill$\Box$

\begin{lemma}\label{cha:Y}
For any  $\delta\in(0,-\frac{\lambda_0^*}{\lambda_0})$, there exists a constant $C>0$,
such that for $|\theta|$ sufficiently large,
\begin{equation*}\|\psi(i\theta)\|_\infty\le C|\theta|^{-\delta}.\end{equation*}
\end{lemma}

\noindent \textbf{Proof:}
For any $\delta\in(0,-\frac{\lambda_0^*}{\lambda_0})$, there exists $\epsilon\in(0,1)$
such that
\begin{equation}\label{domi-lambda*}(1+\epsilon)e^{\lambda_0^* }\le e^{-\lambda_0\delta}.\end{equation}
It follows from Lemma \ref{barVto0} and Lemma \ref{lem:11} that there exists $j\ge 1$
such that for all  $k\ge j$,
\begin{equation}\label{11.9}
  \sup_{\theta\in[1,e^{\lambda_0}]}\|\overline{V}_k(|\psi(i\theta)|)\|_\infty\le \epsilon
   e^{-\|b\|_\infty}.
\end{equation}
Thus by \ref{lemV} we get that,  for all $\theta\in[1,e^{\lambda_0}]$ and $n\ge 1$,
\begin{align*}
  |\overline{V}_{n+j}(\psi(i\theta))(x)|\le& \overline{V}_{n+j}(|\psi(i\theta)|)(x)=\overline{V}_1\overline{V}_{n+j-1}(|\psi(i\theta)|)(x)\\
  \le& (1+\epsilon)T^*_1(\overline{V}_{n+j-1}(|\psi(i\theta)|))(x).
\end{align*}
By iteration and \eqref{T*} we get that,  for $\theta\in[1,e^{\lambda_0}]$,
we can use \eqref{domi-lambda*} to get that
\begin{align*}
\displaystyle |\overline{V}_{n+j}(\psi(i\theta))(x)|\le& \displaystyle(1+\epsilon)^nT^*_n(\overline{V}_{j}(|\psi(i\theta)|))(x)\\
 \le&\displaystyle(1+\epsilon)^nT^*_n(1)(x)\le (1+c)(1+\epsilon)^n e^{\lambda_0^* n} \langle 1,\psi_0^*\rangle_m\|\phi_0^*\|_\infty\\
 \le&\displaystyle(1+c)\langle 1,\psi_0^*\rangle_m\|\phi_0^*\|_\infty e^{\lambda_0\delta(j+1)} e^{-\lambda_0\delta (n+j)}\theta^{-\delta}.
\end{align*}
Since $\psi(i\theta e^{\lambda_0(n+j)})(x)=\overline{V}_{n+j}(\psi(i\theta))(x)$, we have
\begin{equation*}\|\psi(i\theta)\|_\infty\le C |\theta|^{-\delta}, \qquad \theta\ge e^{\lambda _0 j}. \end{equation*}
Using  $\psi(-i\theta)(x)=\overline{\psi(i\theta)(x)}$,  we get
\begin{equation*}\|\psi(i\theta)\|_\infty\le C |\theta|^{-\delta}, \qquad \theta\le -e^{\lambda_0 j}. \end{equation*}
Summarizing, we get the conclusion of the lemma.
\hfill$\Box$

\begin{prop}\label{prop:density}
For any non-zero $\mu\in {\cal M}_F(E)$, under $\P_{\mu}$,  $Y$ is an absolutely continuous random variable, that is, it has a
density function $g_{\mu}(y)$.
\end{prop}
\noindent\textbf{Proof:}
By Remark \ref{rem-den},
it suffices to prove that the conclusion holds when $\mu=\delta_x$.
It follows from  $\psi(\theta,x)=\overline{V}_t(\psi(\theta e^{-\lambda_0t}))(x)$ that
\begin{equation*}Y=^d e^{-\lambda_0t}\sum_{u\in \mathcal{L}_t}Y^u,\end{equation*}
where $\mathcal{L}_t$ is the collection of particles of the branching Markov process
$Z$ which are alive at time $t$.
Given $Z_t$,
$\{Y^u, u\in \mathcal{L}_t\}$ is a family of independent random varaibles with
$Y^u=^d (Y,\P_{\delta_{\xi_t(u)}})$.

Take $\delta\in(0,-\frac{\lambda_0^*}{\lambda_0})$  and $K>0$ such that $K\delta>1$.
For any Lebesgue null set
$B\subset (0,\infty)$, we have
\begin{equation*}
\P_{\delta_x}(Y\in B)\le P_{\delta_x}(\|Z_t\|\le K)+\sum_{n=K+1}^\infty P_{\delta_x}(\|Z_t\|=n, e^{-\lambda_0t}\sum_{u\in \mathcal{L}_t}Y^u\in B).
\end{equation*}
Given $Z_t$ and $\|Z_t\|=n> K$,
for $|\theta|$ sufficiently large, we have
\begin{equation*}
\left|P_{\delta_x}\big(e^{i\theta \sum_{u\in \mathcal{L}_t}Y^u}|Z_t\big)\right|\textbf{1}_{\|Z_t\|=n}\le C^n|\theta|^{-\delta n},
\end{equation*}
implying that the characteristic function of $\sum_{u\in \mathcal{L}_t}Y^u$ is $L^1$ integrable.
Thus $\sum_{u\in \mathcal{L}_t}Y^u$ has a density function, and hence
\begin{equation*}P_{\delta_x}(\|Z_t\|=n, e^{-\lambda_0t}\sum_{u\in \mathcal{L}_t}Y^u\in B)=0.\end{equation*}
Summarizing the above, we have
\begin{equation*}\P_{\delta_x}(Y\in B)\le P_{\delta_x}(\|Z_t\|\le K).\end{equation*}
Letting  $t\to\infty$, we immediately get
$\P_{\delta_x}(Y\in B)=0$,
that is,  the distribution of  $Y$ is absolutely continuous with respect to the  Lebesgue
measure,  and thus has a density function.
\hfill$\Box$

\begin{prop}\label{prop:posit}
For any non-zero $\mu\in\mathcal{M}_F(E)$, under $\P_{\mu}$, the density function of  $Y$
is strictly positive on $(0,\infty)$.
\end{prop}

\noindent\textbf{Proof:}
Note that $\{Y,\P_{\delta_x}\}$ and  $\{aW^Z, P_{\delta_x}\}$ have the same distribution, where $a>0$ is a constant. By Remark \ref{rem-den}, it suffices to show that,
under  $P_{\delta_x}$, the density function of $W^Z$ is strictly positive on $(0,\infty)$.

It has been proven in \cite[Proposition 5.6]{Hering} that, for branching Markov processes satisfying certain conditions, the density function of $W^Z$ is strictly positive on $(0,\infty)$.
For the branching Markov process $\{Z_t\}$ of this paper, we can use the same argument to
show that the same conclusion holds. We omit the details.
\hfill$\Box$

\smallskip

\noindent\textbf{Proof of Theorem \ref{main-dens}:}
Combining Lemma \ref{density2}, Proposition \ref{prop:density} and Proposition \ref{prop:posit}, we immediately get that,  under $\P_{\mu}$,  the distribution of  $W$
is absolutely continuous on $(0,\infty)$ with density function $f_{\mu}$ satisfying that, for all $y>0$,
\begin{equation*}
  f_{\mu}(y)\ge g_{\mu}(y)\langle v,\mu\rangle e^{-\langle v,\mu\rangle}>0.
\end{equation*}
\hfill$\Box$

\subsection{Proof of Theorem \ref{main-tail}}\label{3.3}

Recall that
\begin{equation*}\epsilon_0=\frac{-\lambda_0^*}{\lambda_0}.\end{equation*}

\noindent\textbf{Proof of Theorem \ref{main-tail}}
 First, we deal with the small value probability problem.\\
It follows from Lemma \ref{barV4} that
\begin{equation*}e^{-\lambda_0^* t}\psi(e^{\lambda_0t},x)= e^{-\lambda_0^* t}\overline{V}_t(\psi(1))(x)\to A(\psi(1))\phi_0^*(x),\end{equation*}
that is,
\begin{equation*}\lim_{\theta\to\infty}\theta^{\epsilon_0}\psi(\theta,x)= A(\psi(1))\phi_0^*(x).\end{equation*}
Simple calculations give that, as $\theta\to\infty$,
\begin{align*}
\P_{\mu}(e^{-\theta W}|W>0)=&
\frac{1}{1-e^{-\langle v,\,\mu\rangle}}\Big(e^{-\langle \Phi(\theta),\mu\rangle}-e^{-\langle v,\,\mu\rangle}\Big)
\\=&\frac{1}{e^{\langle v,\,\mu\rangle}-1}\left(e^{\langle \psi(\theta)v,\mu\rangle}-1\right)
\\\sim&\frac{1}{e^{\langle v,\,\mu\rangle}-1} \langle\psi(\theta)v,\,\mu\rangle.
\end{align*}
Summarizing the above, we get
\begin{equation*}
\lim_{\theta\to\infty}\theta^{\epsilon_0}\P_{\mu}(e^{-\theta W}|W>0)
=\frac{1}{e^{\langle v,\,\mu\rangle}-1}\lim_{\theta\to\infty}\theta^{\epsilon_0}\langle\psi(\theta)v,\,\mu\rangle
=\frac{1}{e^{\langle v,\,\mu\rangle}-1} A(\psi(1))\langle v\phi_0^*,\,\mu\rangle.
\end{equation*}
It follows from the Tauberian theorem that
\begin{equation*}
\lim_{r\to 0}r^{-\epsilon_0}\P_{\mu}(W\le r|W>0)=\frac{1}{e^{\langle v,\,\mu\rangle}-1} A(\psi(1))\langle v\phi_0^*,\,\mu\rangle/\Gamma(\epsilon_0+1).
\end{equation*}
Thus
\begin{equation*}
\lim_{r\to 0}r^{-\epsilon_0}\P_{\mu}(0<W\le r)
=e^{-\langle v,\,\mu\rangle}A(\psi(1))\langle v\phi_0^*,\,\mu\rangle/\Gamma(\epsilon_0+1).
\end{equation*}

Now we deal with the tail probability problem.
Let
\begin{equation*}G(s):=\int_0^s \P_{\mu}(W>r)\,dr.\end{equation*}
Then the Laplace transform of $G$ is
\begin{align*}
  \int_0^\infty e^{-\theta r}\,dG(r) =& \int_0^\infty e^{-\theta r}\P_{\mu}(W>r)\,dr \\
  =&\theta^{-1}\Big(1- \theta\int_0^\infty e^{-\theta r}\P_{\mu}(W\le r)\,dr\Big)\\
  =&\theta^{-1}\Big(1- \P_{\mu}(e^{-\theta W})\Big)\\
  =&\theta^{-1}(1- e^{-\langle \Phi(\theta),\mu\rangle}).
\end{align*}
It follows from \eqref{limit-Phi} that
\begin{equation*}
\lim_{\theta\to0}\theta^{-1}\widetilde{L}(\theta^{-1})^{-1}\Phi(\theta,x)=\lim_{t\to\infty}e^{\lambda_0t}L(t)^{-1}\Phi(e^{-\lambda_0t},x)=\phi_0(x),
\end{equation*}
where $L(t)$ are $\widetilde{L}$ are defined in \eqref{L(t)} and \eqref{tilde-L}.
Hence,
\begin{equation*}
\lim_{\theta\to0}\widetilde{L}(\theta^{-1})^{-1}\int_0^\infty e^{-\theta r}\,dG(r)
=\lim_{\theta\to0}\theta^{-1}\widetilde{L}(\theta^{-1})^{-1}\langle \Phi(\theta),\mu\rangle=\langle \phi_0,\mu\rangle.
\end{equation*}
It follows from the Tauberian theorem that
\begin{equation*}
\lim_{r\to\infty}\widetilde{L}(r)^{-1}G(r)=\langle \phi_0,\mu\rangle.
\end{equation*}
Therefore, by \cite{sen2}, we have
\begin{equation*}
\lim_{r\to\infty}r\widetilde{L}(r)^{-1}\P_{\mu}(W>r)=0.
\end{equation*}
\hfill$\Box$

\begin{singlespace}

\end{singlespace}

\vskip 0.2truein
\vskip 0.2truein

\noindent{\bf Yan-Xia Ren:} LMAM School of Mathematical Sciences \& Center for
Statistical Science, Peking
University,  Beijing, 100871, P.R. China. Email: {\texttt
yxren@math.pku.edu.cn}

\smallskip
\noindent {\bf Renming Song:} Department of Mathematics,
University of Illinois,
Urbana, IL 61801, U.S.A
and School of Mathematical Sciences, Nankai University, Tianjin 300071, P. R. China.
Email: {\texttt rsong@illinois.edu}

\smallskip

\noindent{\bf Rui Zhang:} School of Mathematical Sciences, Capital Normal
University,  Beijing, 100048, P.R. China. Email: {\texttt
zhangrui27@cnu.edu.cn}

\end{doublespace}

\end{document}